\documentclass[a4,12pt]{article}

\usepackage{amssymb,amsmath}
\usepackage{xypic}

\newtheorem{defn}{Definition}[section]
\newtheorem{proposition}[defn]{Proposition}
\newtheorem{corollary}[defn]{Corollary}
\newtheorem{exm}[defn]{Example}
\newtheorem{lemma}[defn]{Lemma}
\newtheorem{theorem}[defn]{Theorem}

\newtheorem{xproof}{{\it Proof. }}

\newenvironment{definition}{\begin{defn}\em}{\end{defn}}
\newenvironment{example}{\begin{exm}\em}{\end{exm}}
\newenvironment{proof}{\begin{xproof}\em}{\end{xproof}}

\def\qed{\hspace{0.3cm}{\rule{1ex}{2ex}}}
\newcommand\ie{i.e.}
\newcommand\eg{e.g.}
\newcommand\cf{cf.}
\newcommand\st{\mid}
\newcommand\then{\&}
\newcommand\Spec{{\rm S}}
\newcommand\complex{\mathbb{C}}
\newcommand\rs{\mathrm{R}}
\newcommand\ls{{\rm L}}
\newcommand\ts{{\rm I}}
\newcommand\spat{{\rm Spat}}
\newcommand\Q{{\cal Q}}
\newcommand{\B}{{\cal B}}
\renewcommand{\2}{{\bf 2}}
\newcommand{\Max}{{\rm Max}\,}
\renewcommand\L{{\cal L}}
\newcommand\V{\bigvee}
\newcommand\points{{\rm pt}}
\newcommand\opens{{\cal O}}
\newcommand\pidl{{\downarrow}}
\newcommand\opp[1]{{#1}^{\rm op}}
\newcommand\SpMax{{\rm SpMax}\,}
\newcommand\wi{\eqslantless}
\newcommand\colimit{\lim\limits_{\textrm{\raisebox{.6ex}{$\longrightarrow$}}}}
\newcommand\locales{{\bf Loc}}
\newcommand\ucstar{{\bf C}^{\bf *}}
\newcommand\cucstar{{\bf CC}^{\bf *}}
\newcommand\kreglocales{{\bf KRLoc}}
\newcommand\quantales{{\bf Qu}}
\newcommand\uquantales{{\bf Qu}_e}
\newcommand\squantales{{\bf Qu}_1}
\newcommand\stsquantales{{\bf Qu}_2}
\newcommand\gquantales{{\bf GelQu}}

\begin{document}

\title{On quantales and spectra of C*-algebras\thanks{
The first author was supported by the Grant Agency of the Czech
Republic under grant 201/99/0310. The second
author was supported by a grant
from the Natural Sciences and Engineering Research Council of Canada. The third author
was supported by the Funda\c{c}\~{a}o para a Ci\^{e}ncia e Tecnologia through the
Research Units Funding Program and through grant POCTI/1999/MAT/33018
``Quantales". The fourth author was supported by the Ministry of Education
of the Czech Republic under project MSM 143100009.}}
\author{David Kruml$^1$, Joan Wick Pelletier$^2$, Pedro Resende$^3$, 
and Ji\v{r}\'{\i} Rosick\'y$^1$
\vspace*{2mm}\\
\small\it $^1$Department of Algebra and Geometry, Faculty of Sciences, 
Masaryk University,\vspace*{-2mm}\\
\small\it Jan\'{a}\v{c}kovo n\'{a}m.\ 2a, 66295 Brno, Czech Republic\\
\small\it $^2$Department of Mathematics and Statistics, York University,\vspace*{-2mm}\\
\small\it North York, Ontario, M3J 1P3, Canada\\
\small\it $^3$Departamento de Matem{\'a}tica, Instituto Superior T{\'e}cnico,
\vspace*{-2mm}\\
\small\it Av. Rovisco Pais, 1049-001 Lisboa, Portugal}

\date{~}

\maketitle

    \vspace*{-1,3cm}
\begin{abstract}
We study properties of the quantale spectrum $\Max A$ of an arbitrary
unital C*-algebra $A$. In particular we show that
the spatialization of $\Max A$ with respect to one of the notions
of spatiality in the
literature yields the locale
of closed ideals of $A$ when
$A$ is commutative. We study under general conditions functors with 
this
property, in addition requiring that colimits be preserved, and we conclude 
in this case that the spectrum of $A$ necessarily coincides with
the locale of closed 
ideals of the commutative reflection of $A$. Finally, we address
functorial properties
of $\Max\!$, namely studying (non-)preservation of limits and colimits. 
Although $\Max\!$ is not an equivalence of categories, 
therefore not providing a direct generalization of Gelfand duality to the 
noncommutative case, it is a faithful
complete invariant of unital C*-algebras.
\vspace{0.1cm}
\vspace{0.1cm}
\end{abstract}

\section{Introduction}\label{introduction}

Gelfand duality tells us that any unital commutative C*-algebra is, up to isomorphism, an algebra of
continuous functions on a compact Hausdorff space. More than that, it gives us a dual equivalence
between the category of compact Hausdorff spaces and the category of unital commutative C*-algebras.
Whereas the first fact is at the basis of seeing noncommutative C*-algebras as if they were
algebras of continuous functions on ``noncommutative spaces",
as in noncommutative
geometry~\cite{NCG}, the latter provides a notion of what a category of such
spaces should be,
suggesting that a ``continuous map" of noncommutative spaces $f:A\rightarrow B$ may be defined to be
a $*$-homomorphism of C*-algebras
$f:B\rightarrow A$. Of course, the
generalization of the classical concepts to the noncommutative world can often be done in more
than one way, and for noncommutative spaces there is evidence, \eg\ stemming from
the notion of strong Morita equivalence~\cite{Rie}, that $*$-homomorphisms are not the
right notion of continuous map, or at least not the only one, for commutative C*-algebras which
are strongly Morita equivalent are always isomorphic, but in general strong Morita
equivalence is coarser than isomorphism. For instance, a noncommutative C*-algebra may
be strongly Morita equivalent to a commutative one, indicating that in spite of being
noncommutative it corresponds to a classical (\ie, commutative) space. See~\cite[Ch.\ II.2]{NCG}
for examples.

In the present paper we will focus on another way of studying maps between
noncommutative spaces, based on the observation that any Hausdorff space
$X$ can be recovered up to homeomorphism from its lattice of open sets $\opens(X)$. Concretely, $X$ is
homeomorphic to the space $\points(\opens(X))$ whose points are the maximal proper open sets, \ie, the sets
of the form
$X-\{x\}$ for each point $x\in X$, and whose opens are the sets of the form
$V_U=\{P\in\points(\opens(X))\st U\not\subseteq P\}$, for each open set $U\in\opens(X)$.
This is an instance of the general construction of a space from a locale (see \eg\
\cite{Joh82}), and together with the well known fact that
$\opens(X)$ is order-isomorphic to the lattice $\ts{(C(X))}$ of closed ideals of $C(X)$ this
tells us that the spectrum of a unital commutative C*-algebra $A$ can be obtained directly from
its lattice of closed ideals $\ts(A)$, and indeed the space
$\points(\ts(A))$ is precisely the usual space of maximal ideals of $A$.
Pursuing this line of thought, we can
reformulate Gelfand duality in terms of an equivalence of categories between the category
of unital commutative C*-algebras and the category of compact regular locales (which
is dually equivalent to the category of compact Hausdorff spaces --- see~\cite{Joh82} and
references therein), or, if we wish to work in an arbitrary Grothendieck topos, an
equivalence between the category of unital commutative C*-algebras and the category of
compact completely regular locales~\cite{BanMul1,BanMul2}. This suggests that in the
noncommutative case such a lattice-theoretic approach may provide both an alternative
notion of noncommutative space and an alternative notion of map between noncommutative
spaces.

In the case of a noncommutative C*-algebra
$A$, there are three immediate generalizations of the locale of closed ideals, namely
the
lattices $\ts(A)$, $\rs(A)$, and $\ls(A)$, respectively of closed two-sided, right-sided, or
left-sided ideals. One of them, $\ts(A)$, is in fact isomorphic to the lattice of open
sets of the primitive spectrum of $A$ (the Jacobson topology), and therefore cannot capture that
topological information in $A$ which is of a noncommutative nature. Concerning the other two, it
is irrelevant whether we consider right- or left-sided ideals from here 
on, and we will
address only the right-sided ones. The lattice
$\rs(A)$ is in general not distributive, but the multiplication of closed right ideals obtained
by taking the topological closure of the usual product of ideals,
\[I\then J=\overline{IJ}=\overline{\{a_1 b_1+\ldots + a_n b_n \st a_i\in 
I,\ b_i\in J\}}\;,\]
is distributive over joins:
\[I\then\V_j J_j =\V_j I\then J_j\;,\hspace*{1,5cm}\biggl(\V_i 
I_i\biggr)\then J=\V_i I_i\then J\;.\]
Besides, in the commutative case we have $I\then J=I\cap J$, so the multiplicative structure of
$\rs(A)$ is a generalization of the meets of the locales of closed ideals of the commutative
case. In~\cite{Mul86} it was in fact proposed that the resulting notion of {\em quantale\/},
consisting of a complete lattice $Q$ equipped with an associative multiplication $\then$
that satisfies the two distributivity laws above, might provide a suitable algebraic
framework within which to study the spectrum of C*-algebras, without
implying that such quantales should necessarily be those of closed
right ideals (although some papers seem to have taken this point of
view~\cite{BorBos89,BorBos86,BorRosBos89,Cahiers,RosSpat}). Indeed,
in~\cite{BorRosBos89} it was shown that any postliminal C*-algebra can be recovered from
its quantale $\rs(A)$, but the same is not true for an arbitrary unital C*-algebra.
In~\cite{Curacao} Mulvey has proposed that the spectrum of a unital C*-algebra
$A$ should be taken to be the quantale of all closed subspaces of
$A$, which he called
$\Max A$, whose multiplication is defined in an analogous way to that of $\rs(A)$
(see~\S\ref{background}).

Various properties of the quantales $\Max A$ have meanwhile been
studied~\cite{MulPel1,MulPel2,PasRos}, and the
``continuous maps" between quantales are thought of as (the opposites of)
quantale
homomorphisms
$h:Q\rightarrow Q'$,
\ie, functions which preserve joins and the multiplicative structure,
\[h(a\then b)=h(a)\then h(b)\;,\hspace*{1,5cm}h\biggl(\V_i a_i\biggr)=\V_i h(a_i)\;,\]
and satisfying other requirements (``strength" or ``unitality") meant to force quantale
homomorphisms between locales to be precisely the homomorphisms of locales
(see~\S\ref{background}). This leads to some good candidates for categories of noncommutative
spaces in the quantalic sense, and the main purpose of our present paper is to address, albeit
preliminarily, the way in which these categories relate to a ``standard" category of
noncommutative spaces, namely the category of unital C*-algebras, leaving
matters such as Morita
equivalence to future work.

Our paper is based on three main ideas. First, motivated by the fact that in locale theory
only some locales are spatial --- \ie, isomorphic to topologies of actual spaces --- we ask
the question of whether the quantales of the form $\Max A$ are spatial. In particular, in
\S\ref{spatiality} we examine the notion of spatiality studied in~\cite{PelRos,Pas,Kru},
based on a strong embedding into a product of certain simple quantales. This notion
provides an elegant generalization of the notion of spatiality for locales, but
$\Max A$ is not usually spatial in this sense, as an example due to Kruml~\cite{PasRos}
shows. In an attempt to remedy this situation, a notion of ``spatialization"
is introduced and the ``spatial spectrum" $\SpMax A$ is studied. Although its functorial
properties are unsatisfactory, it evolves that for commutative C*-algebras $A$, $\SpMax A$
and the localic spectrum of $A$, \ie\ the locale of closed ideals, coincide.

This last fact suggests our second idea, which is to investigate
the existence of a functor $\Spec$ from C*-algebras to quantales agreeing with the localic
spectrum on commutative algebras. In the hope that $\Spec$ will be a potential equivalence
of categories, producing an extension of Gelfand duality to noncommutative C*-algebras, we
also require that it preserve colimits, a necessary condition for a left adjoint.
But, as we shall see in Theorem~\ref{thm:trivialS}, such a functor $\Spec$ 
just gives us the 
locale of closed ideals of the commutative reflection of $A$, which of 
course is completely uninteresting from the point of view of 
noncommutative topology --- for instance, $\Spec(A)$ is trivial for any simple C*-algebra (\ie, 
an algebra without proper closed two-sided ideals), which includes all
the algebras of matrices $M_{n}(\complex)$ and various interesting examples 
that arise in noncommutative geometry, such as the noncommutative space 
of Penrose tilings of the plane~\cite[Ch.\ II.3]{NCG}.

In~\cite{MulPel2} Mulvey and Pelletier formulated another notion of spatiality, based on
the notion of a discrete von Neumann quantale, motivated by work of Giles and
Kummer~\cite{GK} and Akemann~\cite{Ake70}, in which $\Max A$ is spatial for all C*-algebras
$A$. It turns out that this notion is characterized by a right embedding (\ie\ an embedding
on right-sided elements) into a product of certain simple quantales, 
which partly revives the role
played by right-sided quantales
in earlier papers by
various authors~\cite{BorBos89,BorBos86,BorRosBos89,Cahiers}.
It follows that it is not the spectrum $\Max A$ itself that coincides
with the localic spectrum when $A$ is commutative, but rather its localic
reflection, $\ts(\Max A)$, which coincides. The difficulties met in
~\S\S~\ref{spatiality} and~\ref{extending} when trying to describe a spatial spectrum 
reinforce the idea that the revised notion
of spatiality is necessary --- at least as far as studying C*-algebras is 
concerned, for other applications of quantales may lead to different 
notions of spatiality~\cite{ResJPAA}.

Our acceptance of the spatiality of $\Max A$ leads us
in~\S\ref{functorial} to the third part 
of this paper,
namely, to study the properties of the functor $\Max\!$ itself, and
we attempt to describe the category of
quantales that arises as image of the functor $\Max\!$. Although the immediate
conclusion is that the categories are rather different, since in 
particular the functor $\Max\!$ does not have any adjoints, it is the 
case that $\Max\!$ is faithful. 
Furthermore, it is possible, using results
in~\cite{GK,MulPel1},
to recover $A$ from $\Max A$ for any unital C*-algebra $A$, \ie, 
$\Max\!$ is a complete invariant of unital C*-algebras, and thus
it seems to have promising properties.

\section{Background}\label{background}

We begin by recalling the definitions of quantale theory which have come
to
play an important role in analysis and that of the motivating example,
$\Max A$. By a {\em sup-lattice\/} we shall always mean a complete lattice, and by a {\em
sup-lattice homomorphism\/} we mean a map $f$ between two sup-lattices which preserves arbitrary
joins (but not necessarily meets):
\[f\biggl({\V_i x_i}\biggr)=\V_i f(x_i)\;.\]

\begin{definition}
By a {\em quantale\/} $Q$ is meant a sup-lattice
together with an
associative product, $\then$,  satisfying
\[a \then \biggl(\V_i b_{i}\biggr)\ = \V_i (a \then b_{i}) \]
and
\[\biggl(\V_i a_{i}\biggr) \then b = \V_i (a_{i} \then b) \]
for all $a, b, a_{i}, b_{i} \in Q$. We always denote the bottom element $\V\emptyset$ of $Q$ by
$0$. The quantale
$Q$ is said to be {\em unital\/} provided
that there exists an element $e \in Q$ for which
\[e \then a = a = a \then e \]
for all $a \in Q$. The top element of $Q$ is denoted by $1$, and a quantale whose top
satisfies
$1\then 1=1$ is {\em strong\/}. A quantale $Q$ for which $0=1$ is {\em trivial\/}.
\end{definition}

\begin{example}\label{exm:QS}
Let $S$ be a sup-lattice. Then the set $\Q(S)$ of all the sup-lattice endomorphisms of $S$ is a
unital quantale whose joins are calculated pointwise, $(\V_i f_i)(x)=\V_i f_i(x)$, whose multiplication is
composition,
$f\then g=g\circ f$, and whose unit is the identity map ${\rm id}_S$.
Such a quantale $\Q(S)$ is an example of a simple quantale in the sense
of~\cite{Pas}.
\end{example}

\begin{definition}
An element $a \in Q$ of a quantale $Q$ is said to be {\em right-sided\/} if
\[a \then 1 \leq a \]
for all $a \in Q$.
Similarly, $a \in Q$ is said to be {\em left-sided\/} if
\[1 \then a \leq a \]
for all $a \in Q$, and {\em two-sided\/} if it is both right- and left-sided. The sets of right-,
left-, and two-sided elements of
$Q$ are respectively denoted by $\rs(Q)$, $\ls(Q)$, and $\ts(Q)$.
\end{definition}

We remark that for a unital quantale $Q$ an element $a\in Q$ is right-sided 
if and only if $a=a\then 1$,
and thus $\rs(Q)= Q\then 1$. Similarly, $\ls(Q)=1\then Q$ and 
$\ts(Q)=1\then Q\then 1$ if $Q$ is unital.

\begin{definition}
By an {\em involutive quantale\/} is meant a quantale $Q$ together with an
involution, $^*$,
satisfying the conditions that
\[a^{**} = a\;, \]
\[(a \then b)^* = b^* \then a^*\;, \mbox{ and } \]
\[\biggl(\V_i a_{i}\biggr)^* = \V_i a_{i}^* \]
for all $a, b, a_{i} \in Q$.

By a {\em Gelfand quantale\/} is meant a quantale
$Q$ which is unital, involutive, and which
satisfies the condition
\[a \then a^* \then a = a \]
for each right-sided (equivalently, left-sided) element $a \in Q$.
\end{definition}

Motivated by Mulvey's consideration of the closed right ideals $\rs (A)$ of a
C*-algebra $A$ as a quantale, Rosick\'{y}~\cite{Cahiers} introduced a new
multiplication and resulting structure (called quantum frame) on $\rs(A)$ in an
attempt to find a complete invariant of $A$. While the quantum frame
$\rs(A)$
failed to provide such an invariant and failed too to be functorial in
the general
sense, the importance of the new multiplication and of the role of the
closed right
ideals in the representation of  C*-algebras in turn led Mulvey~\cite{Curacao} to
propose the quantale $\Max A$,
the lattice of all
closed linear subspaces of $A$, as the spectrum of $A$. Specifically,

\begin{definition}
For a noncommutative C*-algebra $A$ with identity, by the {\em
spectrum\/} of $A$ is meant the unital involutive quantale $\Max A$ of all closed
linear subspaces
of $A$ together with the product, $\then$, defined by setting
\[M \then N = \overline{M N} = \overline{\{a_1 b_1+\ldots + a_n b_n \st 
a_i\in M,\ b_i\in N\}}\;,\]
to be the closure of the product of linear subspaces for each $M, N \in
\Max
A$, and the join, $\V$, defined by taking
\[ \V_{i} M_{i} = \overline{\sum_{i} M_{i}} \]
to be the closure of the sum of linear subspaces for each family $M_{i}
\in
\Max A$. The involution is calculated pointwise,
$M^{*}=\{a^{*}\st a\in M\}$,
and the unit of $\Max A$ is given by the closed linear subspace
generated
by the identity of $A$.
\end{definition}

We note that $\Max A$ is actually a Gelfand quantale~\cite{Curacao}.

\begin{definition}
By a {\em quantale homomorphism\/} $\varphi:Q\rightarrow Q'$ is
meant a sup-lattice homomorphism
which also preserves the operation of product, $\then$, of
the
quantales. The map is {\em unital\/} provided that
$e' = \varphi(e)$,
where $e$ and $e'$ are the respective units, {\em pre-unital\/} if $e'\le
f(e)$,\footnote{In~\cite{MulPel1} the unital homomorphisms are what we are now calling
pre-unital.} and {\em strong\/} if
$1' = \varphi (1)$,
 where $1$ and $1'$ are the respective top
elements of $Q$ and $Q'$. A  homomorphism of involutive quantales is
{\em involutive\/} if
it preserves the involution.

We denote various categories of quantales as follows:
$\quantales$ is the category of quantales and quantale homomorphisms; by adding a subscript
``$e$", as in $\uquantales$, we restrict to unital quantales and unital homomorphisms; by adding
a subscript ``$1$", as in $\squantales$, we restrict to strong quantales and strong homomorphisms;
by adding a superscript ``$*$", as in $\quantales^*$, we restrict to involutive quantales and
involutive homomorphisms. (For instance, $\uquantales^*$ is the category of unital involutive
quantales with unital and involutive homomorphisms.)
\end{definition}

\begin{example}
The mapping $A\mapsto \Max A$ extends to a
functor from the category $\ucstar$ of unital $C^*$-algebras to the category
$\uquantales^*$ (or to its full
subcategory
$\gquantales$ of Gelfand quantales)~\cite{MulPel1}: if $f:A
\rightarrow B$, then $\Max f:\Max A \rightarrow \Max B$ assigns to each $M \in A$ the
closure of
$f[M]=\{f(a)\in B\st a\in M\}$.
\end{example}

\begin{definition}
Let $Q$ be a quantale and $S$ a sup-lattice. A {\em representation\/} of $Q$ on $S$ is a homomorphism
$r:Q\rightarrow\Q(S)$. This makes $S$ a {\em right module\/} over $Q$ in an obvious sense, and we write
$x\cdot a$ for the action of $a\in Q$ on $x\in S$, \ie, $x\cdot a=r(a)(x)$.

The representation $r$
is said to be {\em irreducible\/} if the only principal ideals $I\subseteq S$ which are
invariant under the representation are $\{0_S\}$ and $S$ itself; that is, if the
condition $I\cdot Q\subseteq I$ implies either
$I=\{0_S\}$ or $I=S$.
\end{definition}

\begin{proposition}\label{prop:irreducibility}
\begin{enumerate}
\item A representation of $Q$ on $S$ is irreducible if and only if
$x\cdot 1_Q\le x$ implies that $x=0_S$ or $x=1_S$,
for all
$x\in S$.
\item Any strong representation is irreducible.
\item In the case of unital quantales, a pre-unital representation
is irreducible if and only if it is
strong.
\end{enumerate}
\end{proposition}

\begin{proof}
1. Obvious from the fact that the condition $x\cdot 1_Q\le x$ is equivalent to $\pidl
x\cdot Q\subseteq \pidl x$ [since ${\pidl x}\cdot Q=\pidl(x\cdot 1_Q)$], $x=0_S$ is
equivalent to
$\pidl x=\{0_S\}$, and
$x= 1_S$ is equivalent to $\pidl x=S$.\\
2. Assume that a representation of a quantale $Q$ on
$S$ is strong. In module notation this means that
$x\cdot 1_Q=1_S$ for all $x\neq 0$ in $S$. Now assume $x\cdot 1_Q\le x$ for some $x\in
S$.
If $x\neq 0$ we must have $x=1_S$, and thus the representation is
irreducible.\\
3. Assume that a representation of a unital quantale $Q$ on $S$ is
pre-unital and irreducible. We have
$x\cdot 1_Q\cdot 1_Q\le x\cdot 1_Q$, and thus $x\cdot 1_Q$ is either $0$ or $1_S$,
by irreducibility. But
$x\cdot 1_Q\ge x\cdot e_Q\ge x$ due to pre-unitality, and thus $x\cdot 1_Q=0$
if and only if $x=0$. Hence,
$x\cdot 1_Q=1_S$ for all $x\neq 0$.
\qed
\end{proof}

\begin{example}\label{exm:Hilbertreps}
Let $A$ be a unital C*-algebra and
$\pi:A\rightarrow\B(H)$
a representation
of $A$ on a Hilbert space
$H$. Then~\cite{MulPel1} a representation $\tilde\pi$ of $\Max A$
on $\L(H)$ can be
obtained canonically by setting, for each closed subspace $V\in\Max A$
and each closed subspace
$W\in\L(H)$,
\[\tilde\pi(V)(W)=\overline{\langle\{\pi(a)(x)\st a\in V,\ x\in 
W\}\rangle}\;,\]
where $\langle X\rangle$ denotes the linear span of a subset 
$X\subseteq H$. If $\pi$ is irreducible then so is 
$\tilde\pi$~\cite{MulPel1}. We shall refer to such an irreducible 
representation as a {\em Hilbert point\/} of $\Max A$.
\end{example}

To conclude this section we recall that a {\em locale\/}
is a sup-lattice $L$
satisfying the distributivity law
\[a\wedge\V_i b_i=\V_i a\wedge b_i\]
for all $a,b_i\in L$, and that a {\em locale
homomorphism\/} is a sup-lattice homomorphism $f:L\rightarrow L'$ that also preserves the top and
the binary meets (equivalently, the finite meets):
\begin{eqnarray*}
f(a\wedge b)&=& f(a)\wedge f(b)\;,\\
f(1)&=& 1'\;.
\end{eqnarray*}
Hence, any locale is a unital involutive quantale whose multiplication is meet, $\wedge$, whose
involution is trivial, $a^*=a$, and whose unit is the top, $1$. The category of locales and their
homomorphisms is denoted by $\locales$ (beware that in~\cite{Joh82} 
$\locales$ stands for our $\opp\locales$).

\section{Spatialization}\label{spatiality}

In this section we address the extent to which quantales $\Max A$ can be viewed as
``spaces" in the sense of being strongly embedded into a product of certain simple
quantales, a notion introduced in~\cite{PelRos} for involutive quantales and further studied
in~\cite{Pas,Kru} for arbitrary quantales. Examples based on 
Theorem~\ref{thm:kruml} will show
that $\Max A$ is usually not spatial in this sense. Furthermore, we 
will see that attempts to
make the spatialization construction functorial lead to undesirable consequences.

A natural definition of spatial involutive quantale is the following, meant to
generalize that of spatial locale:

\begin{definition}\label{def:spat}
Let $Q$ be a quantale. We say that $Q$ is {\em spatial\/} with respect to a family of
homomorphisms (``points") $\{p_i:Q\rightarrow Q_i\}_i$ if, for all $a,b\in Q$, $a=b$ whenever $p_i(a)=p_i(b)$ for
all $i$.
\end{definition}

Of course, this notion of spatiality depends upon the class of points considered.
In~\cite{PelRos} a definition was given for an involutive quantale $Q$, taking as ``points"
the involutive irreducible representations $p:Q\rightarrow\Q(S)$ on a sup-lattice $S$ with
a duality (\ie, a sup-lattice $S$ equipped with an order 2 antitone automorphism), where $\Q(S)$ is an involutive quantale due to the 
duality of $S$~\cite{MulPel0}, and by involutive representation is meant a 
representation which is an involutive homomorphism. In other words, $Q$ was said to be spatial if it could be strongly embedded into
a product of simple involutive quantales of the form $\Q(S)$. In~\cite{MulPel1}
the notion of ``point" of a unital involutive quantale 
was further refined, although no notion of spatiality was offered
in that article. The object of~\cite{MulPel1} was to ensure when $Q$ is of the form $\Max
A$ that points would correspond precisely to irreducible representations of the C*-algebra
$A$.

A practical characterization of the notion of spatiality above was given in~\cite{Kru} by
means of prime elements of a quantale, where the points of a quantale $Q$ can be
taken to be the irreducible representations
$p:Q\rightarrow \Q(S)$. This characterization will be useful for us so we recall it here.

\begin{definition}\label{def:prime}
Let $Q$ be a quantale. An element
$p\in Q$ is {\em prime\/} if the condition $a\then 1 \then b\le p$ implies $a\le p$ or $b\le p$,
for all
$a,b\in Q$.
\end{definition}

\begin{theorem}[Kruml]\label{thm:kruml}
Let $Q$ be a quantale. Then $Q$ is spatial (with respect to irreducible representations) if and
only if each
$a\in Q$ is a meet of prime elements of $Q$.
\end{theorem}

With this characterization it is easy to see that a quantale of the form $\Max A$ is not
necessarily spatial in the sense of Definition~\ref{def:spat}, even if we take the points
to be all the irreducible representations:

\begin{example}\label{exm:M2Cnotspatial}
The ``one-point" C*-algebra of $2\times 2$ matrices $A=M_2(\complex )$ has a
non-spatial spectrum. To see this let
\begin{eqnarray*}
    P&=&\{ (^a_b\ ^b_c)\mid 
a,b,c\in\complex\}\;,\\
R&=&\{
(^a_0\ ^b_0)\mid a,b\in\complex\}\;,\\
L&=&\{ (^a_b\ ^0_0)\mid a,b\in\complex\}\;.
\end{eqnarray*}
We obtain $R\then A\then L=\{ (^a_0\
^0_0)\mid a\in\complex )\}\subseteq P$ but $R,L\not\subseteq P$. Hence, $P$ is maximal but
not prime, and thus it cannot be obtained as a meet of primes, which means $\Max
A$ is not spatial.
\end{example}

In order to understand better the (non-)spatiality of $\Max A$, let us study the
case when $A$ is commutative. First we define a notion of spatialization for (not
necessarily commutative) quantales that is analogous to that of locales, although it turns
out, as we shall see, not to possess the same good functorial properties as in the locale
case.

\begin{definition}
Let $\{p_i\}_i$ be a family of points of a quantale $Q$. The {\em spatialization\/}
of $Q$ with respect to $\{p_i\}_i$ is the quotient $Q/{\sim}$, where $a\sim b$ if and
only if $p_i(a)=p_i(b)$ for all $i$. We denote by $\spat_Q:Q\rightarrow Q/{\sim}$ the natural
quotient homomorphism.
\end{definition}

\begin{theorem}\label{thm:spat}
Let $A$ be a commutative unital C*-algebra. Then the following quantales coincide:
\begin{enumerate}
\item the spatialization of $\Max A$ with respect to the Hilbert 
points (\cf\ Example~\ref{exm:Hilbertreps});
\item the spatialization of $\Max A$ with respect to all the irreducible representations;
\item the locale of closed ideals $\ts (\Max A)$ (\ie, the ``classical" localic spectrum of
$A$).
\end{enumerate}
Moreover, the natural map $\spat_{\Max A}$ is given by $a\mapsto 1 \then 
a\then 1$.
\end{theorem}

\begin{proof}
We have a quotient $\Max A\rightarrow\ts(\Max A)$ that sends each
$a\in\Max A$ to $1 \then a\then 1$ (the localic reflection of $\Max A$).
Since the points in 1 form a proper subclass of the points in 2,
the quantale
in 1 is a quotient of that in 2, and it suffices to show that the
quantale in 2 is a quotient of $\ts (\Max A)$.

Let us represent
$A$ as the algebra of continuous functions $C(X)$ on a compact Hausdorff
topological space $X$. Consider a prime element $P$ in $\Max A$ and let
$u,v\in X$ be arbitrary points of $X$. They have disjoint open neighborhoods,
say $U$ and $V$. Let $I_F=\{f\in C(X)\mid(\forall x\in F)(f(x)=0)\}$ be
the closed ideal corresponding to a closed set $F$. Then $(X-U)\cup (X-V)=X$
implies $I_{X-U}I_{X-V}=\{0\}\subseteq P$. Since $P$ is prime, we have that
$I_{X-U}\subseteq P$ or $I_{X-V}\subseteq P$. Since the points $u,v$ have been
chosen arbitrarily, we deduce that there is at most one point $x_0$ which has no
open neighborhood $W$ such that $I_{X-W}\subseteq P$. For the other elements
$x\in X-\{ x_0\}$ we have a system of open sets $U_x$ (not containing $x_0$) such
that $I_{X-U_x}\subseteq P$ for every $x$, hence $P$ contains the maximal
ideal $I_{\{ x_0\}}=\bigcup_x I_{X-U_x}$.

Now, $I_{\{x_0\}}$ is maximal not only as an ideal but also as a subspace of $C(X)$, and thus
$P=I_{\{x_0\}}$ because
$P$ is prime and thus $P\neq C(X)$. In order to verify that $I_{\{x_0\}}$ is a maximal subspace
assume that
$P$ contains a function
$f$ such that
$f(x_0)\neq 0$, and let $g\in C(X)$. Define a
function
$h\in C(X)$ by
\[h(x)= g(x)-\frac{f(x)}{f(x_0)} g(x_0)\;.\]
Then $h(x_0)=0$, and thus $h\in P$, whence $g\in P$ because $P$ is a linear space and $g$ is the
linear combination
$h+\lambda f$, with $\lambda=g(x_0)/f(x_0)$.
Since $g$ is an arbitrary function in
$C(X)$ it follows that $P=C(X)$.

We conclude that every prime element of
$\Max A$ is two-sided and maximal, and thus the spatialization of $\Max A$ is
a
locale~\cite{PasRos} and coincides with $\ts(\Max A)$. Suppose now that $P$ is any closed
subspace of $A$. Then the image in the spatialization corresponds to the
least closed ideal
containing $P$, which equals $A\then P\then A$. \qed
\end{proof}

This theorem provides us with a large supply of C*-algebras with a non-spatial spectrum, and
indeed
$\complex^2$ is probably the simplest example of such an algebra ---
$\Max\complex^2$ is infinite, whereas its spatialization $\ts(\Max\complex^2)$ is a four element
Boolean algebra; alternatively, one can see immediately that
$\{ (x,x)\mid x\in\complex\}$ is maximal but not prime, whence $\Max\complex^2$ is not spatial
by the same argument of Example~\ref{exm:M2Cnotspatial}.

Now we turn to some functorial properties of spatialization. First, for each quantale
$Q$ let
$Q/{\sim}$ be the spatialization with respect to a class of points 
${\cal C}_{Q}$ of 
which we require only that the two following properties hold, in the  
case when 
$Q=\Max A$ for some unital C*-algebra $A$:
\begin{enumerate}
    \item ${\cal C}_{\Max A}$ contains all the Hilbert points of $\Max 
    A$
    (\cf\ Example~\ref{exm:Hilbertreps}).
    \item If $A$ is commutative then $(\Max A)/{\sim}=\ts(\Max A)$.
\end{enumerate}
For each unital C*-algebra $A$
let us write
$\SpMax A$ for the quotient
$(\Max A)/{\sim}$. It is natural to ask whether $\SpMax\!$ can be extended to a functor
from C*-algebras to spatial quantales, as this would yield an interesting spectrum because it
would coincide with the classical localic spectrum of the commutative case. In fact the answer
turns out to be negative, at least if the family
$\{\spat_{\Max A}\}_A$ is required to be a natural transformation from
$\Max\!$ to
$\SpMax\!$. This means that spatialization of quantales in this sense is not
as well behaved
as spatialization for locales, a not surprising fact due to the mismatch between
arbitrary quantale homomorphisms, which are not necessarily strong, and irreducible
representations, which are. The following results provide a clear-cut argument.

\begin{lemma}\label{lemma:pushout}
Let $f:\complex^2 \rightarrow M_2(\complex)$ be the embedding of unital C*-algebras
given
by $f(a,b)=(^a_0\ ^0_b)$. Consider the pushout $Q$ of $\Max f$ and
$\spat_{\Max\complex^2}$ in
$\quantales$:
\[\xymatrix{ \Max\complex^2 \ar[d]_{\spat_{\Max\complex^2}} \ar[rr]^{\Max f}
&& \Max M_2(\complex) \ar[d]^p \\ \ts(\complex^2) \ar[rr]_q && Q\;.}
\]
Then $Q$ is trivial.
\end{lemma}

\begin{proof}
Let $E,F$ be the one-dimensional subspaces of $\complex^2$ spanned by
$(1,1)$ and $(1,-1)$, respectively. Then
$\spat_{\Max\complex^2}(E)=\spat_{\Max\complex^2}(F)=\complex^2$. Since the subspaces spanned by the
matrices $(^1_0\ ^0_1)$ and
$(^1_{0}\ ^{~0}_{-1})$
are the images $(\Max f)(E)$ and $(\Max f)(F)$,
they must be equalized by $p$. Let
\begin{eqnarray*}
L &=& \left\{{\begin{matrix}a & a \\ b & b\end{matrix}\mid
a,b\in\complex}\right\},\\ R &=& \left\{{\begin{matrix}c & d \\ c &
d\end{matrix}\mid c,d\in\complex}\right\}.
\end{eqnarray*}
Then from 
\[\begin{matrix} a & a \\ b & b\end{matrix}\begin{matrix} 1 & 0 \\ 0 &
1\end{matrix}\begin{matrix} c & d \\ c & d\end{matrix} = \begin{matrix} a &
a \\ b & b \end{matrix}\begin{matrix} c & d \\ c & d \end{matrix} =
2\begin{matrix} ac & ad
\\ bc & bd\end{matrix}\;,\]
\[\begin{matrix} a & a \\ b & b\end{matrix}\begin{matrix} 1 & 0 \\ 0 &
-1\end{matrix}\begin{matrix} c & d \\ c & d\end{matrix} = \begin{matrix} a &
-a \\ b & -b \end{matrix}\begin{matrix} c & d \\ c & d \end{matrix} =
\begin{matrix} 0 & 0
\\ 0 & 0 \end{matrix}\;,\]
it follows that $L\then (\Max f)(E)\then R=M_2(\complex)$, while $L\then (\Max f)(F)\then R=\{
0\}$. Hence,
$p(M_2(\complex))=p(\{ 0\})$ and, since the diagram commutes,
$q(\complex^2)=q(\{0\})$. This means that both $p$ and $q$ have the trivial image
$\{0\}\subseteq Q$, and thus the trivial quantale $Q=\{0\}$ is the pushout in $\quantales$.
\qed
\end{proof}

\begin{lemma}\label{lemma:spmaxM2Cnottrivial}
$\SpMax M_2(\complex)$ is not trivial.
\end{lemma}

\begin{proof}
In order to show this it suffices to find a non-zero point of $\Max M_2(\complex)$. The
isomorphism $\imath:M_2(\complex)\rightarrow \B(\complex^2)$ is an irreducible representation of
$M_2(\complex)$ on $\complex^2$ and gives rise to a Hilbert point
$\tilde \imath:\Max M_2(\complex)\rightarrow\Q(\L(\complex^2))$, as in
Example~\ref{exm:Hilbertreps}. This point is clearly non-zero, as for instance the
identity matrix
spans a subspace of $M_2(\complex)$ that acts as the identity on
$\L(\complex^2)$.
\qed
\end{proof}

\begin{theorem}
There is no functorial extension of $\SpMax\!$ such that the family $\{\spat_Q\}_Q$ is a natural
transformation.
\end{theorem}

\begin{proof}
Assume that there is a functorial extension of $\SpMax\!$ such that $\{\spat_Q\}_Q$ is a natural
transformation. Then we should have a
commuting diagram
\[\xymatrix{ \Max\complex^2 \ar[d]_{\spat_{\Max\complex^2}} \ar[rr]^{\Max f}
&& \Max M_2(\complex) \ar[d]^{\spat_{\Max M_2(\complex)}} \\ \SpMax\complex^2 \ar[rr]_{\SpMax f}
&&
\SpMax M_2(\complex)\;,}
\]
where $f$ is defined as in
Lemma~\ref{lemma:pushout} and $\SpMax\complex^2=\ts(\complex^2)$.
It follows that $\spat_{\Max M_2(\complex)}$ must factor through the 
pushout of Lemma~\ref{lemma:pushout}, which is $0$. Thus, $\SpMax M_2(\complex)$ must be
trivial itself since $\spat_{\Max M_2(\complex)}$ is onto, which leads to a contradiction
with Lemma~\ref{lemma:spmaxM2Cnottrivial}.
\qed
\end{proof}

\section{Extending the localic spectrum}\label{extending}

The previous section showed that for quantales of the form $\Max A$, with $A$
commutative, the spatialization with respect to ``reasonable'' points is just
the locale of closed ideals of $A$. This
suggests the possibility that a functor from C*-algebras to ``noncommutative spaces" should
have the property of coinciding with the localic spectrum in the commutative case.
In the previous section we have seen that there are difficulties in trying to define such a
functor from $\Max\!$ via spatialization, so we now address this idea in a different way.
We shall define a functor $\Spec$ from unital C*-algebras to quantales, requiring only that
it should coincide with the localic spectrum in the commutative case and that it preserve
colimits, which has the advantage of opening the way to finding an adjunction between C*-algebras
and quantales. We will see, however, that under very general assumptions the spectrum of a
C*-algebra $A$ given by such a functor is necessarily a locale, in fact 
coinciding with the localic spectrum of the commutative reflection of 
$A$, which renders it uninteresting from the point of view of 
noncommutative topology, as was stated in the introduction.

From here on let $\Spec$ denote a functor from $\ucstar$ to some category of involutive
quantales
${\cal Q}$ that contains the category of locales, but which otherwise will remain
unspecified for the moment. An important
observation is that any unital C*-algebra is a colimit of commutative unital
C*-algebras (more precisely, a quotient of a coproduct of commutative unital
C*-algebras), and thus if
$\Spec$ preserves all the colimits of $\ucstar$ and yields the localic spectrum in the
commutative case, it follows that
the spectrum of an arbitrary unital C*-algebra is a colimit of locales in ${\cal Q}$. In
particular, the initial object of the category should be preserved, which means that the two
element chain $\2=\Spec(\complex)$ should be initial in ${\cal Q}$.
In this section we study the effects of these restrictions when ${\cal Q}$ is
made to coincide with certain subcategories of $\quantales^*$. Notice that
$\quantales^*$ itself is immediately ruled out because $\2$ is not initial in it.
Hence, we focus on two particular subcategories, namely $\uquantales^*$ and
$\squantales^*$, both of which have
$\locales$ as a full subcategory.

First we recall that a {\em strictly two-sided quantale\/} is
a unital quantale whose elements are all two-sided, or, equivalently, a unital quantale
whose unit is the top (see \cite{Ros90}).
We can think of such a quantale as being a ``non-idempotent locale", for
locales are exactly the strictly two-sided quantales whose multiplication is idempotent.
We denote by
$\stsquantales$ the category of strictly two-sided quantales and strong
(=unital) homomorphisms, regarded as involutive quantales with the 
trivial involution $a^*=a$. This is a full subcategory both of $\uquantales^*$ 
and $\squantales^*$.

\begin{lemma}\label{lem:2scolims}
    $\stsquantales$ is closed under 
    colimits, both in $\uquantales^*$ and in
    $\squantales^*$.
\end{lemma}

\begin{proof}
    $\stsquantales$ has the same initial object of both categories, 
    it is obviously closed under coequalizers (the quotient of a 
    strictly two-sided quantale is strictly two-sided), so it remains to verify 
    that it is closed under coproducts.
    Let $\{Q_{i}\}_{i}$ be a family of strictly two-sided quantales.
    We first consider the coproduct in $\uquantales^*$.
    The unit of each quantale $Q_{i}$ is preserved by the respective coprojection, 
    and for each $i$ we have $a\le e$ for all $a\in Q_{i}$, whence it 
    follows that all the generators of the coproduct are below the unit.
    Hence, the unit of the coproduct is the top because in any unital
    involutive quantale the set of elements below the unit is closed under multiplication,
    joins and the involution.
Now let us consider $\squantales^*$.
For each $i$ the top of $Q_{i}$ is preserved by the respective
coprojection, and thus the new top is a
neutral element for any finite non-empty multiplication of elements from the
quantales $Q_{i}$. Also, we are not considering multiplications with zero length because we
do not require our coproduct quantale to be unital, so all the elements of the
coproduct are either finite non-empty multiplications or joins
of these (the involution cannot generate new elements because it is 
preserved by the coprojections), and thus the top is a neutral element. \qed
\end{proof}

From here on we write $\mathcal Q$ to mean either $\uquantales^*$ or 
$\squantales^*$ whenever the distinction is not relevant.

\begin{corollary}
    A colimit of locales in $\mathcal Q$ is a strictly 
    two-sided quantale.
\end{corollary}

Now recall that a locale $L$ is \emph{regular\/} if each of its elements 
$a\in L$ is the join of those $a'\in L$ which are \emph{well inside\/} 
$a$ (written $a'\wi a$),
\[a=\V\{a'\in L\st a'\wi a\}\;,\]
where $a'\wi a$ means that $a'\wedge b=0$ and $a\vee b= 1$ for some 
$b\in L$.

\begin{lemma}
    Let $\{L_{i}\}_{i}$ be a family of locales, with all but at most 
    one of the $L_{i}$ regular. Then the
    coproduct $\coprod_{i}L_{i}$ in $\mathcal Q$ is calculated in $\locales$.
\end{lemma}

\begin{proof} If the family is empty the coproduct is the initial 
object, which is in $\locales$. Otherwise, we already know that $L=\coprod_{i}L_{i}$ is a strictly two-sided
    quantale, so it is
enough to show that it is also idempotent.
First, we remark that all the generators of $L$ are idempotent because all
the coprojections $\iota_{i}:L_{i}\to L$ are quantale homomorphisms.
Furthermore, a join of idempotents is
idempotent, and if two idempotents $a$ and $b$ commute then $a\then b$ is
idempotent: 
\[a\then b\then a\then b=a\then a\then b\then b=a\then b\;.\] 
Hence, to finish the proof we will show that $\iota_{i}(a)$ and $\iota_{j}(b)$ commute for 
all $a\in
L_{i}$ and $b\in L_{j}$. By assumption at least one of $L_{i}$ and 
$L_{j}$ is regular, so assume that $b$ is a join $\V_{k}b_{k}$
with $b_{k}\wi b$ for all $k$. For each $k$ let then $c_{k}\in L_{j}$
be an element such that $b_k\wedge c_k=0$ and $b\vee c_k= 1$.
We obtain
\begin{align*}
\iota_{i}(a)\then \iota_{j}(b_k) &= \iota_{j}(b\vee c_k)\then \iota_{i}(a)\then \iota_{j}(b_k) \\
              &= \iota_{j}(b)\then \iota_{i}(a)\then \iota_{j}(b_k)\vee \iota_{j}(c_k)\then
	      \iota_{i}(a)\then \iota_{j}(b_k) \\
              &\leq \iota_{j}(b)\then \iota_{i}(a)\vee \iota_{j}(c_k\wedge 
              b_k) = \iota_{j}(b)\then 
              \iota_{i}(a)\;,
\end{align*}
and thus $\iota_{i}(a)\then \iota_{j}(b)=\bigvee \iota_{i}(a)\then 
\iota_{j}(b_k)\leq \iota_{j}(b)\then \iota_{i}(a)$. 
Similarly, multiplying by
$\iota_{j}(b\vee c_k)$ on the right we obtain
$\iota_{j}(b)\then \iota_{i}(a)\leq \iota_{i}(a)\then \iota_{j}(b)$.
\qed
\end{proof}

\begin{corollary}
    The category $\kreglocales$ of compact regular locales is closed under colimits in 
    $\mathcal Q$.
\end{corollary}

\begin{proof}
    It is well known that
    $\kreglocales$ is closed under colimits in $\locales$ (see, \eg, 
    \cite{Joh82}),
    which in turn is obviously closed under 
    coequalizers in $\mathcal Q$. The previous lemma shows that the 
    coproduct of compact regular locales in $\mathcal Q$ is also
    calculated in $\locales$ 
    and, therefore, in $\kreglocales$. \qed
\end{proof}

\begin{theorem}\label{thm:trivialS}
Let $\Spec:\ucstar\rightarrow{\cal Q}$ be a functor that preserves colimits, and such
that $\Spec(B)=\ts(B)$ for every commutative unital C*-algebra $B$. Let also
$A$ be an arbitrary unital C*-algebra. Then $\Spec(A)$ is the locale of closed ideals
of the commutative reflection of $A$.
\end{theorem}

\begin{proof}
    Let $A$ be a unital C*-algebra, obtained as a colimit
    $\colimit B_{i}$, in $\ucstar$, of commutative unital C*-algebras.
Then $\Spec(A)$ is a colimit $\colimit \ts(B_{i})$ of compact regular
locales in $\mathcal Q$, and by the above results this colimit is 
calculated in $\kreglocales$. The commutative reflection 
functor
of unital C*-algebras preserves colimits because it is a left 
adjoint, and thus the commutative reflection of $A$ can be 
obtained as a ``similar'' colimit as $A$, but calculated in 
$\cucstar$ instead of $\ucstar$. Furthermore, the localic spectrum 
$\ts:\cucstar\rightarrow\kreglocales$ preserves colimits because it 
is an equivalence of categories, and thus the localic spectrum of the 
commutative reflection of $A$ equals $\colimit 
\ts(B_{i})=\Spec(A)$. \qed
\end{proof}

\section{Functorial properties of $\Max\!$}\label{functorial}

In recent work~\cite{MulPel2} a different approach is taken to spatiality of involutive
quantales, essentially motivated by~\cite{GK}. In particular it is shown that
any quantale $\Max A$ is spatial with respect to
the Hilbert points, provided one changes the definition of spatiality so that it
applies to right-sided elements only, following from the idea~\cite{Ake70,GK} that the
topology of $\Max A$ is provided by the right-sided elements. Explicitly, a quantale
$Q$ is spatial in this sense if and only if for any distinct right-sided elements $a$
and $b$ there is a point
$p$ of
$Q$ such that
$p(a)\neq p(b)$. Of course, a general notion of spatialization in
this sense is likely to suffer
from the same problems we met in \S\ref{spatiality}, but
this need not concern us
because $\Max A$ is already spatial in the new sense.
Instead we may focus on the functor $\Max\!$ itself.
We claim that $\Max\!$ has some pleasing properties,
namely that it is faithful and a complete
invariant of unital C*-algebras. Hence, it is an interesting candidate
via which to search for an equivalence of categories that
generalizes the Gelfand duality of the commutative case,
although this is certainly far from being straightforward because as we
will see $\Max\!$ does not have any adjoints.

Let us first see that $\Max\!$ does not preserve colimits, and thus it does
not have a right
adjoint:

\begin{theorem}\label{thm:sumsnotpres}
$\Max:\ucstar\rightarrow\uquantales$ does not preserve coproducts.
\end{theorem}

\begin{proof}
Let $\complex^2+\complex^2$
denote the sum of $\complex^2$ with itself in the category of unital C*-algebras. The sum
in the category of commutative unital C*-algebras is $\complex^4$ and there is a quotient
morphism
$\complex^2+\complex^2\rightarrow\complex^4$. Hence, we have an onto unital homomorphism of
quantales $\Max(\complex^2+\complex^2)\rightarrow\Max\complex^4$ (notice
that $\Max f$ is onto whenever $f$ is). If the coproduct
$\Max\complex^2+\Max\complex^2$ in $\uquantales$ were isomorphic to
$\Max(\complex^2+\complex^2)$ then there would be, by composition, an onto homomorphism
$h:\Max\complex^2+\Max\complex^2\rightarrow\Max\complex^4$ that would make the
following diagram commute,
\[
\xymatrix{\Max\complex^2\ar[drrr]^{\Max\gamma_1}\ar[d]_{\kappa_1}\\
\Max\complex^2+\Max\complex^2\ar[rrr]^h&&&\Max\complex^4\\
\Max\complex^2\ar[urrr]_{\Max\gamma_2}\ar[u]^{\kappa_2}}
\]
where $\kappa_1$ and $\kappa_2$ are the coprojections of the coproduct in $\uquantales$,
and $\gamma_1$ and $\gamma_2$ are the coprojections of the coproduct in the category of
commutative unital C*-algebras. In order to see that sums are not preserved we will
show that $h$ is not onto. This is equivalent to showing that the unital involutive subquantale $Q$ of
$\Max\complex^4$ generated by the union of the images of $\Max\gamma_1$ and
$\Max\gamma_2$ is not the whole quantale $\Max\complex^4$. For this we look at
$\gamma_1$ and $\gamma_2$ in detail. These coprojections are the duals of the projections
of the product of two discrete topological spaces with two points, $\{x,y\}\times\{x,y\}$.
We denote the complex valued functions $f:\{x,y\}\rightarrow\complex$ by
ordered pairs $(f(x),f(y))$. Similarly, we denote the functions on the product space,
$f:\{(x,x),(x,y),(y,x),(y,y)\}\rightarrow\complex$, by 4-tuples
$(f(x,x),f(x,y),f(y,x),f(y,y))$. The two projections of the product thus have the following
duals:
\begin{eqnarray*}
(z,w)&\stackrel{\gamma_1}\mapsto&(z,z,w,w)\;,\\
(z,w)&\stackrel{\gamma_2}\mapsto&(z,w,z,w)\;.
\end{eqnarray*}
In any quantale $\Max A$ we have, for all $a,b\in A$, $\langle a\rangle\then\langle
b\rangle=
\langle ab\rangle$, where $\langle a\rangle$ and $\langle b\rangle$ are the one dimensional
subspaces spanned by $a$ and $b$. The one dimensional subspaces are the atoms of $\Max
A$, and an atom $\langle a\rangle$ cannot be obtained as the join of different atoms, for
if $a\neq b$ then $\langle a\rangle\vee\langle b\rangle$ is a two dimensional subspace.
Consider now the one dimensional subspace $\langle(1,0,1,1)\rangle$ of $\Max\complex^4$.
For this subspace to be in $Q$ it is necessary that
$(z,z,w,w)(z',w',z',w')=(1,0,1,1)$ for some $z,w,z',w'\in\complex$ (the involution is of no use, of course).
But then
either
$z$ or
$w'$ must be $0$, whence $(zz',zw',wz',ww')$ must contain at least two zeros, and we
conclude that $h$ is not onto. \qed
\end{proof}

It is worth remarking that $\Max\!$ does not preserve limits, either, so it does not
have a left adjoint:

\begin{theorem}
$\Max\!$ does not preserve products.
\end{theorem}

\begin{proof}
$\Max\complex$ is the two element chain $\2$, and thus
$\Max\complex\times\Max\complex$ is the four element boolean algebra. However,
$\complex^2$
contains infinitely many subspaces, whence
$\Max(\complex\times\complex)\not\cong\Max\complex\times\Max\complex$.
\qed
\end{proof}

\begin{theorem}\label{thm:faithfulMax}
$\Max\!$ is faithful.
\end{theorem}

\begin{proof}
This is a corollary of Proposition~2.3 in~\cite{Cahiers}, for if
$f:A\rightarrow B$ is a morphism of unital C*-algebras the map $\rs f: \rs A\rightarrow\rs
B$ on closed right ideals can be defined by $\rs f(J)=\Max f(J)\then B$. \qed
\end{proof}

Finally,
from the results 
in~\cite{GK,MulPel1} it follows~\cite{ResFields} that $\Max\!$
is a complete invariant 
of unital C*-algebras
(see also~\cite{BorRosBos89,Cahiers,PasRos}):

\begin{theorem}\label{reconstruction}
Let $A$ and $B$ be unital C*-algebras. Then $A$ and $B$ are isomorphic 
if and only if $\Max A$ and $\Max B$ are isomorphic as unital involutive 
quantales.
\end{theorem}

\section*{Acknowledgements}

The authors thank Chris Mulvey for his comments on an earlier draft of 
this paper.

\end{document}